\documentclass[11pt]{amsproc}\usepackage{amsmath}

\usepackage{amsfonts}
\usepackage{epsfig}
\theoremstyle{plain}

\newtheorem{corollary}{Corollary}

\newtheorem{definition}{Definition}
\newtheorem{example}{Example}

\newtheorem{lemma}{Lemma}

\newtheorem{proposition}{Proposition}
\newtheorem{remark}{Remark}

\newtheorem{theorem}{Theorem}
\numberwithin{equation}{section}

\input{tcilatex}

\begin{document}
\title[Newton maps for quintic polynomials]{Newton maps for quintic
polynomials}
\author{Francisco Balibrea}
\address{Facultad de Matem\'{a}ticas, Universidad de Murcia Campus de
Espinardo, 30100 Murcia, Spain }
\email{balibrea@um.es}
\author{J. Orlando Freitas}
\address{Dept. of Mathematics and Engineering, University of{\small \ }%
Madeira, Campus Universit\'{a}rio da Penteada, 9000-390 Funchal, Portugal }
\email{orlando@uma.pt}
\author{J. Sousa Ramos}
\address{Department of Mathematics, Instituto Superior T\'{e}cnico, Av.
Rovisco Pais 1, 1049-001 Lisbon, Portugal }
\email{sramos@math.ist.utl.pt}
\urladdr{http://www.math.ist.utl.pt/\symbol{126}sramos}
\date{}

\begin{abstract}
The purpose of this paper is to study some properties of the Newton maps
associated\ to real quintic polynomials. First using the Tschirnhaus
transformation, we reduce the study of Newton's method for general quintic
polynomials to the case $f(x)=x^{5}-c~x+1$. Then we use symbolic dynamics to
consider this last case and construct a kneading sequences tree for Newton
maps. Finally, we prove that the topological entropy is a monotonic
non-decreasing map with respect to the parameter $c$.
\end{abstract}

\maketitle

\section{Introduction and motivation}

The classical problem of solving equations has substantially influenced the
development of mathematics throughout the centuries and still has several
important applications to the theory and practice of present-day computing.
\textquotedblleft Solving the quintic\textquotedblright\ is one of the few
topics in mathematics which has been of enduring and widespread interest for
centuries.

We cannot solve the general polynomial equation of fifth degree or higher
using radicals. Consequently, methods for estimating numerical solutions of
equations as simple as polynomials are necessary. On a mundane level,
numerical methods can be used to find the zeros (real or complex) to any
required degree of accuracy. This is a useful practical method.

Many mathematical problems can be reduced to compute the solutions of $%
f(x)=0 $, and Newton's method 
\begin{equation}
x_{n+1}=N_{f}(x_{n})=x_{n}-\frac{f(x_{n})}{f^{\prime }(x_{n})},\,\text{\ \ }%
n=0,1,2,...  \label{eq401}
\end{equation}%
is the most common algorithm to solve this problem. The geometric
interpretation of the Newton's method is well known. In such a case $x_{n+1}$
is the point where the tangent line $y-f(x_{n})=f^{\prime }(x_{n})(x-x_{n})$
to the graph of $f(x)$ at the point $\left( x_{n},f(x_{n})\right) $
intersects the $x $-axis.

The fundamental property of Newton's method is that it transforms the
problem of finding roots of $f(x)$ into the problem of finding attracting
fixed points of the associate Newton's map $N_{f}(x)$.

We may ask however for the set of all points $x_{0}$ from which the Newton's
method is converging to a solution.

In 1829, Cauchy \cite{Cauchy29}\ first proved a convergence theorem which
does not assume any existence of a solution. Under standard assumptions, the
Newton's method is locally convergent in a suitable disk centered at the
solution. The possibility that a small change in $x_{0}$ can cause a drastic
change in convergence indicates the nasty nature of the convergence problem.

A key notion in the study of discrete dynamical systems is that of \textit{%
chaos}\ and \textit{sensitive dependence on initial conditions}. There have
been several definitions of chaos, for example Devaney's definition \cite%
{Devaney89}. In one-dimensional case Devaney's definition is equivalent to
the existence of a dense orbit and another criterion is considering the
system as chaotic whether the entropy is strictly positive.

A detailed treatment of the cubic polynomial case can be seen in \cite{R}
and a complete description of its combinatorics is given in \cite{SSR99}.

We study the quintic function $%
f(x)=x^{5}+c_{1}x^{4}+c_{2}x^{3}+c_{3}x^{2}+c_{4}x+c_{5}$.

In the first place we reduce the number of parameters. We now discuss the
necessary algebraic aspects of this reduction. As it is well-known, the
equation 
\begin{equation*}
x^{5}+c_{1}x^{4}+c_{2}x^{3}+c_{3}x^{2}+c_{4}x+c_{5}=0,
\end{equation*}%
with arbitrary coefficients $c_{j}$, can be transformed to the \textit{%
Bring-Jerrard} type 
\begin{equation*}
x^{5}+a~x+b=0,
\end{equation*}%
by a Tschirnhaus transformation \cite[pp. 212-214]{Dickson26}.

Tschirnhaus's transformation reduces the $n^{th}$ degree polynomial equation

\begin{equation*}
c_{0}x^{n}+c_{1}x^{n-1}+...+c_{n-1}x+c_{n}=0
\end{equation*}%
to one with up to three fewer terms 
\begin{equation*}
x^{n}+b_{4}x^{n-4}+...+b_{n-1}x+b_{n}=0
\end{equation*}%
by transforming the root as follows 
\begin{equation*}
x_{j}=\gamma _{4}x_{j}^{4}+\gamma _{3}x_{j}^{3}+\gamma _{2}x_{j}^{2}+\gamma
_{1}x_{j}+\gamma _{0},\text{ \ }(j=1,...,n)
\end{equation*}%
where the $\gamma _{j}$ can be expressed in radicals in terms of the $a_{j}$%
. Thus every quintic can be transformed into one of the form 
\begin{equation}
x^{5}+b_{4}~x+b_{5}=0.  \label{eq411}
\end{equation}%
The $b_{j}$ can ultimately be expressed in radicals in terms of the $a_{j}$ 
\cite{Weisstein99}.

\begin{remark}
The resulting expressions are really complicated. For a general quintic with
symbolic coefficients they require a lot of computation and storage. However
this is not a problem of present-day computer.
\end{remark}

\begin{figure}[th]
\vspace{5mm} \centerline{\epsfig{file=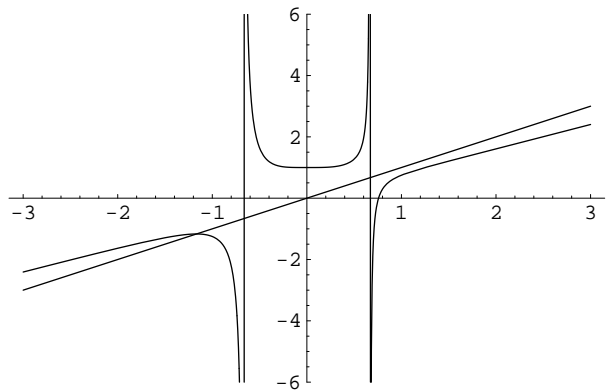,width=8cm}}
\caption{Typical graph of $N_{f_{c}}$ with one fixed point.}
\label{fig401}
\end{figure}

The outline of the paper is as follows. In section 2 we use topological
conjugacy to study the Newton map for quintic polynomials of the form $%
f(x)=x^{5}+a~x+b$. Apparently, the idea of conjugation is important to
understand the iteration of $N_{f}(x).$ Indeed, it was E. Schr\"{o}der who
observed the importance of conjugations, mainly to obtain a convenient form
for calculations, (see Peitgen and Haeseler \cite{Peitgen88}). The concept
of conjugations has proven extremely useful in the modern theory of
iteration and we use it. With this idea we reduce the general case to most
interesting having the form $f_{c}(x)=x^{5}-c~x+1$ whose picture is shown in
Figure \ref{fig401}.

In section 3, using standard symbolic dynamics, we introduce the
admissibility rules of the sequences associated to Newton maps $N_{f}$. Then
we study the structure of the set of admissible sequences. The techniques of
symbolic dynamics are based on the notions of the kneading theory for
one-dimensional multimodal maps, (see Milnor and Thurston \cite{MT86}). We
construct a kneading sequences tree for $N_{f_{c}}$.

In section 4 we are devoted to consider the topological entropy. An
immediate consequence of these results is the exact computation of the
topological entropy which is made in this section. Our main result is the
last theorem of the section.

The connection between kneading theory and subshifts of finite type is shown
by using a commutative diagram derived from the topological configurations
associated with $m$-modal maps, (see Lampreia and Sousa Ramos \cite{LSR93}).

\section{Newton maps for quintics}

We study the polynomial function $f(x)=x^{5}+a~x+b.$ As we said before we
use topological conjugacy. It is well known that $f$ and $g$ are
topologically conjugate provided there is an homeomorphisms $\tau $ such
that $f\circ\tau = \tau\circ g$. In such case for $f^{n}$ and $g^{n}$ we
have the same relationship. So $f$ and $g$ are holomorphically the "same"
dynamical systems. Indeed topological conjugacy is a very efficient concept
to carry over difficult dynamical problems in simpler ones. It plays an
important role in the investigation of the dynamics of general
one-dimensional maps \cite[p. 122]{Collet80}.

\begin{proposition}
Let $g(x)=x^{5}+ax+b^{5}$ where $b\neq 0$ and define $f(x)=x^{5}+c~x+1$
where $c=a/b^{4}$. Then $N_{g}$ and $N_{f}$ are topologically conjugate via
the homeomorphism $\tau (x)=x/b.$
\end{proposition}

\begin{proof}
First we calculate 
\[
N_{g}(x)=x-\dfrac{g(x)}{g^{\prime }(x)}=\dfrac{4x^{5}-b^{5}}{5x^{4}+a}
\]%
and 
\[
N_{f}(x)=x-\dfrac{f(x)}{f^{\prime }(x)}=\dfrac{4x^{5}-1}{5x^{4}+c}.
\]%
We have $\tau \circ N_{g}(x)=\dfrac{4x^{5}-b^{5}}{5bx^{4}+ab}$ and $%
N_{f}\circ \tau (x)=\dfrac{4x^{5}-b^{5}}{5bx^{4}+ab},$ so we have

\[
N_{g}=\tau ^{-1}N_{f}\circ \tau ,
\]%
i.e., $N_{g}$ and $N_{f}$ are topologically conjugate.
\end{proof}

Let us see what happens when $b=0$.

\begin{proposition}
Let $g(x)=x^{5}+a^{4}x$,\ $\tau (x)=x/a$, with $a\neq 0$, and $%
p_{+}(x)=x^{5}+x$. Under such conditions $N_{g}$ and $N_{p_{+}}$ are
topologically conjugate by $\tau $. By other hand, if $g(x)=x^{5}-a^{4}x$
and\ $p_{-}(x)=x^{5}-x$ then $N_{g}$ and $N_{p_{-}}$ are topologically
conjugate by $\tau .$
\end{proposition}

\begin{proof}
We can calculate 
\[
N_{g}(x)=x-\dfrac{g(x)}{g^{\prime }(x)}=\dfrac{4x^{5}}{5x^{4}+a^{2}}
\]
and 
\[
N_{p_{+}}(x)=x-\dfrac{p_{+}(x)}{p_{+}^{\prime }(x)}=\frac{4x^{5}}{5x^{4}+1}.
\]
We have 
\[
\tau \circ N_{g}(x)=\dfrac{4x^{5}}{5ax^{4}+a^{5}}
\]
and 
\[
N_{p_{+}}\circ \tau (x)=\dfrac{4x^{5}}{5ax^{4}+a^{5}},
\]
so

\[
N_{g}=\tau ^{-1}N_{p_{+}}\circ \tau .
\]

It is analogous for the second case and we have 
\[
N_{g}=\tau ^{-1}N_{p_{-}}\circ \tau ,
\]%
i.e., $N_{g}$ and $N_{p_{-}}$ are topologically conjugate.
\end{proof}

We must consider now the case $a=0$ which leaves $f(x)=x^{5}$.

\begin{remark}
Last two propositions imply that either the dynamics of Newton's map for the
quintic $g(x)=x^{5}+a~x+b$ are equivalent to the dynamics of Newton's map
for the polynomial family $f_{c}(x)=x^{5}+c~x+1$ or to $g_{a}(x)=x^{5}+a~x.$
Moreover, the Newton's map for function $g_{a}(x)$ is topologically
conjugate to Newton map for one of the three polynomials: 
\[
p_{-}(x)=x(x^{4}-1),\ \ p_{+}(x)=x(x^{4}+1),\ \ \text{or }p_{0}(x)=x^{5}.
\]
\end{remark}

Therefore the study of Newton map for quintic polynomials is reduced to the
case $f_{c}(x)=x^{5}+c~x+1$. Indeed, as $c\in \mathbb{R}$\ we use $%
f_{c}(x)=x^{5}-c~x+1$ instead of $f_{c}(x)=x^{5}+c~x+1$. In this case we
have $f_{c}^{^{\prime }}(x)=5x^{4}-c$.

\begin{figure}[th]
\vspace{5mm} 
\centerline{
\epsfig{file=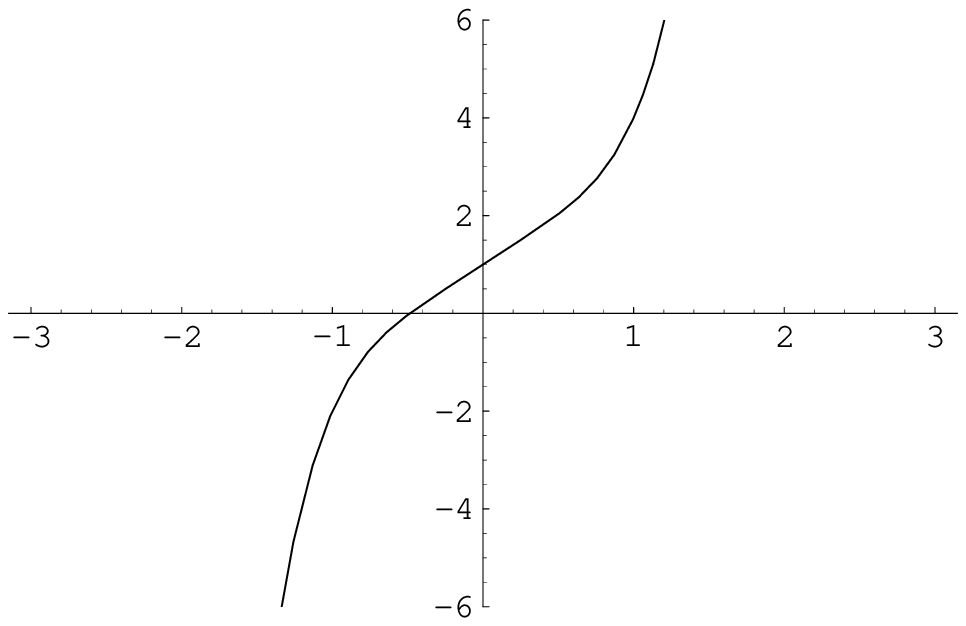,width=6cm}
\epsfig{file=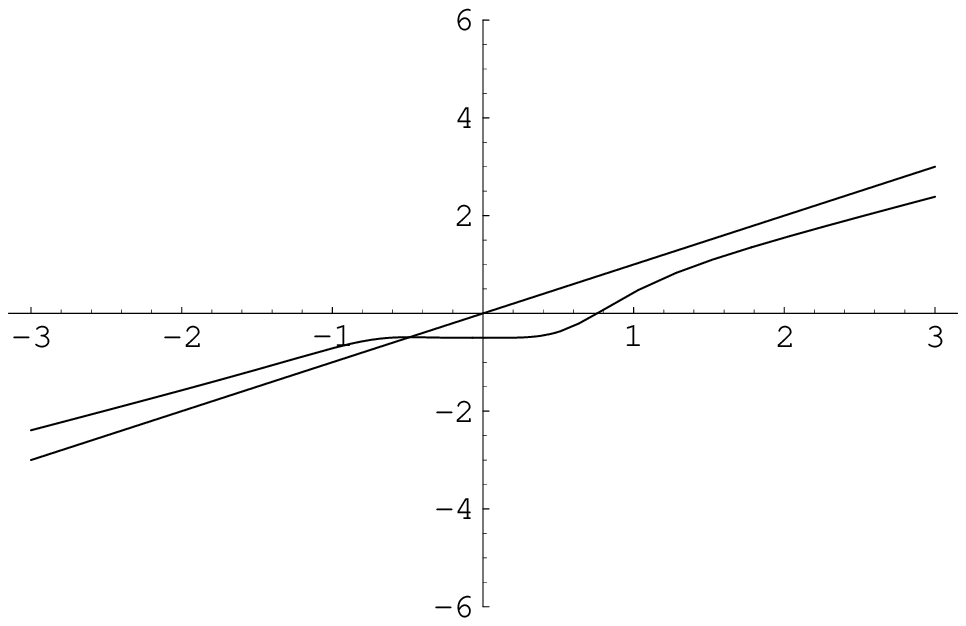,width=6cm}
}
\caption{Map $f_{c}$ and typical map $N_{f_{c}}$ for $c<0$.}
\label{fig402}
\end{figure}

\begin{itemize}
\item When $c<0,$ it is easy to verify that $N_{f_{c}}$ has exactly one real
root and that its \textit{stable} set (the set of points which are forward
asymptotic to it) contains all $\mathbb{R}$ as we see in fig \ref{fig402}.

\item When $c=0$, there is also one real root and its stable set contains
all real numbers except 0.

\item When $c>0$ we have three interesting cases.
\end{itemize}

The polynomial $f(x)$ has a relative maximum at $d_{1}=-\sqrt[4]{c/5}$ and a
relative minimum at $d_{3}=\sqrt[4]{c/5}$.

\begin{figure}[th]
\vspace{5mm} 
\centerline{
\epsfig{file=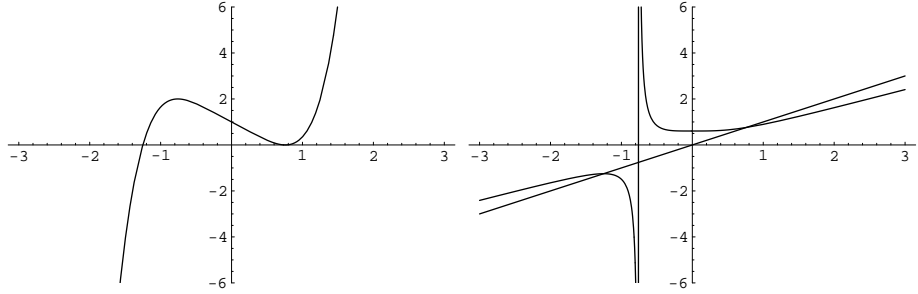,width=6cm}
\epsfig{file=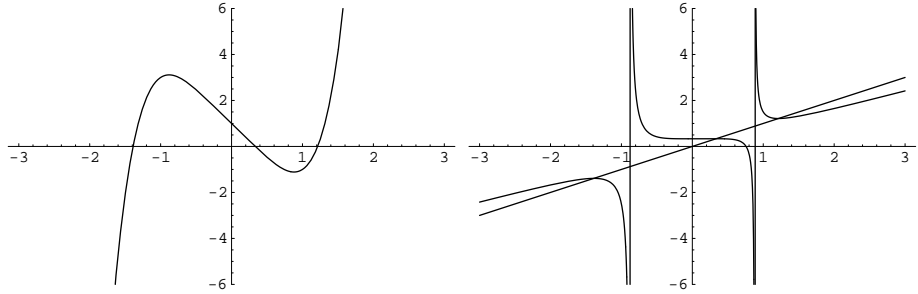,width=6cm}
}
\caption{Map $f_{c}$ and typical map $N_{f_{c}}$ for $c=c_{0}$.}
\label{fig403}
\end{figure}

We note that when $c$ increases, the relatives minimum of $f$ decreases and
the relative maximum increases. When $c=5\times 2^{-8/5}=1.64938$... the
relative minimum is 0, see Figure \ref{fig403}. We denote the parameter $%
c=5\times 2^{-8/5}$ by $c_{0}$.

\begin{figure}[th]
\vspace{5mm} 
\centerline{
\epsfig{file=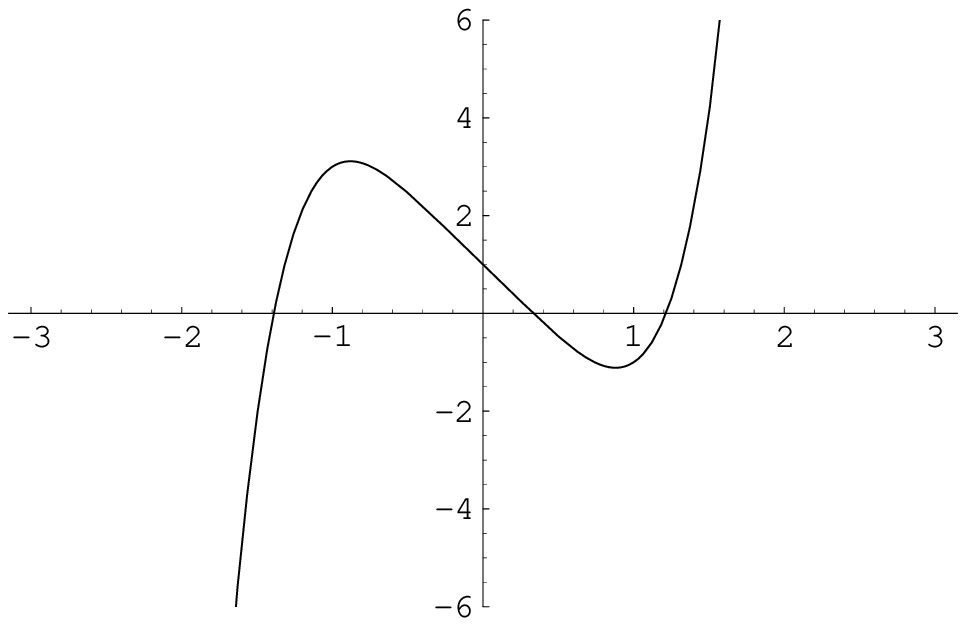,width=6cm}
\epsfig{file=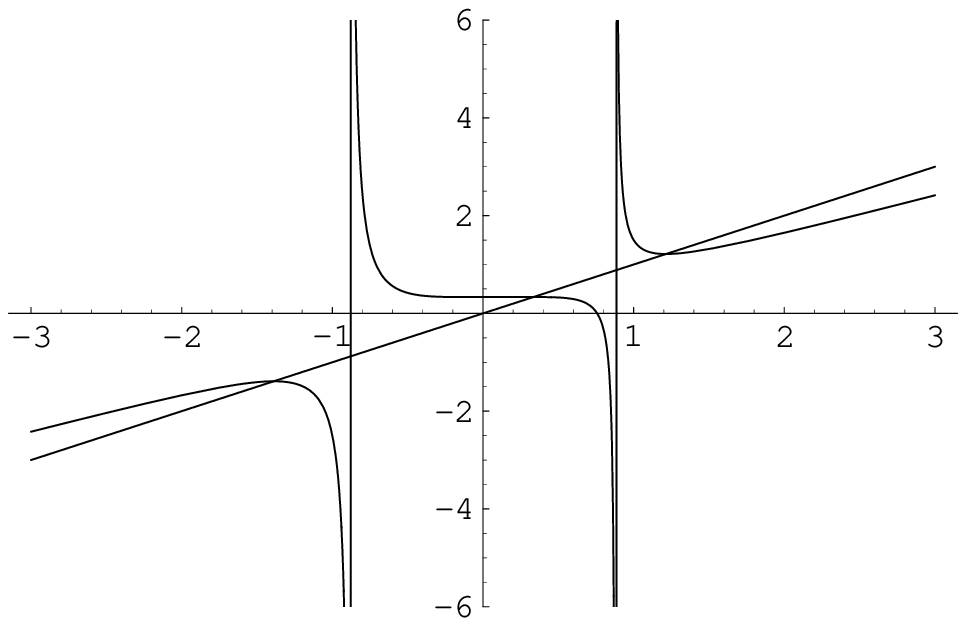,width=6cm}
}
\caption{Map $f_{c}$ and typical map $N_{f_{c}}$ for $c>c_{0}$.}
\label{fig404}
\end{figure}

Note that when $c$ is bigger than $c_{0}$, $f$ has three real roots as
showed in Figure \ref{fig404}.

In last case we can use the following result:

\begin{theorem}[R\'{e}nyi \protect\cite{Renyi50}]
Let $f:\mathbb{R\rightarrow R}$ be defined on $(-\infty ,+\infty )$. Let us
suppose that $f^{\prime \prime }(x)$ is monotone increasing for all $x\in 
\mathbb{R}$ and that $f(x)=0$ has exactly three real roots $a_{i}$ $%
(i=1,2,3) $.

Then the sequence $x_{n+1}=x_{n}-f(x_{n})/f^{\prime }(x_{n})$ converges to
one of the roots for every choice of $x_{0}$ except for $x_{0}$ belonging to
a countable set $E$ of singular points, which can be explicitly given. For
any $\varepsilon >0$ there exists an interval $(x,x+\varepsilon )$
containing three points $y_{i}$ $(x<y_{i}<x+\varepsilon ,$ $i=1,2,3)$ having
the property that if $x_{0}=y_{i} $, then $(x_{n})_{n=0}^{\infty }$
converges to $a_{i}$ $(i=1,2,3)$.
\end{theorem}

The polynomial $f(x)=x^{5}-c~x+1$ has three real roots when $c>c_{0}$ and $%
f_{c}^{^{\prime \prime }}(x)=20x^{3}$ is monotone increasing for all $x\in 
\mathbb{R}$, so we are in the conditions of last theorem.

\begin{figure}[th]
\vspace{5mm} 
\centerline{
\epsfig{file=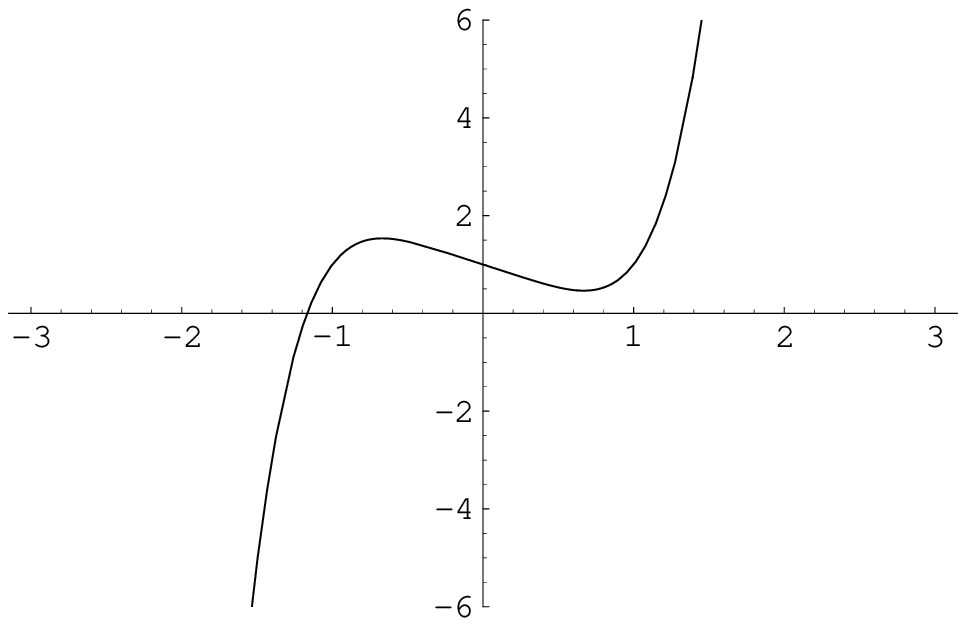,width=6cm}
\epsfig{file=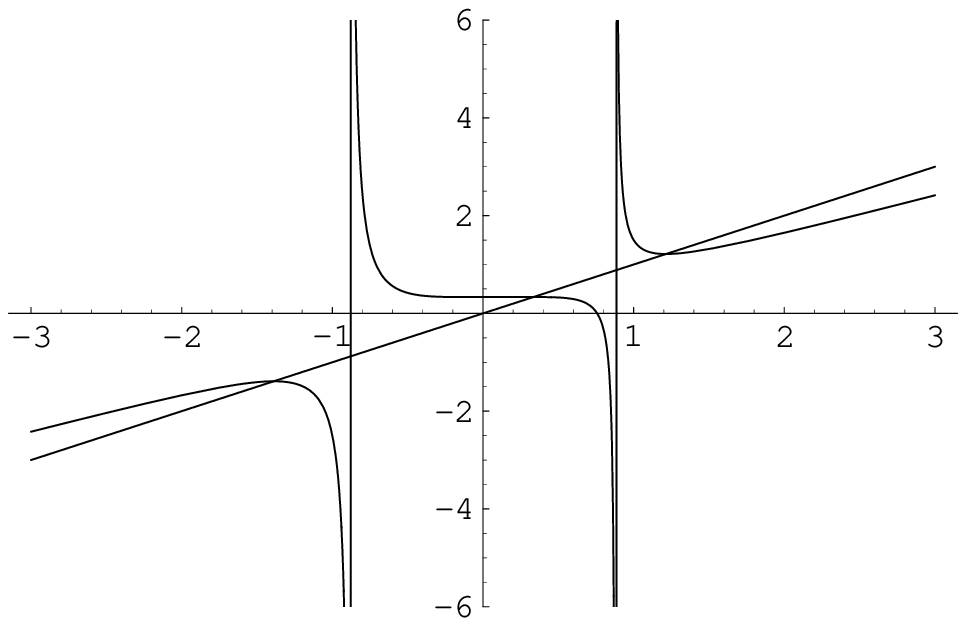,width=6cm}
}
\caption{Map $f_{c}$ and typical map $N_{f_{c}}$ for $0<c<c_{0}$.}
\label{fig405}
\end{figure}

Finally we have the most interesting case, showed in Figure \ref{fig405},
when $c$ is between $0$ and $c_{0}$, in this case $f(x)$ has only one real
root and we denote it by $d_{0}.$

Much of the motivation for the material to be presented comes from the
following theorem due to Fatou \cite{Fatou19}:

\begin{theorem}[Fatou]
Let $R: \overline{\mathbb{C}}\rightarrow \overline{\mathbb{C}}$ be a
rational function with a stable periodic orbit $(R^{n}(z))_{n=0}^\infty$ $%
(|R^{\prime}(z)| < 1)$, then the orbit of at least one critical point $%
\omega $ $(R(\omega) = 0)$ converges to $z$.
\end{theorem}

\begin{proof}
A complete and detailed proof of this fact can be found in \cite{Blanchard84}%
.
\end{proof}

This theorem has an important implication for the family of mappings $%
N_{f_{c}}(x)$, because there is not another stable periodic orbit except
that of the critical point of $f_{c}(x)$.

Now we consider $f_{c}(x)=x^{5}-c~x+1$, so 
\begin{equation*}
N_{f_{c}}^{\prime }(x)=\frac{f_{c}^{^{\prime \prime }}(x)f_{c}(x)}{\left(
f_{c}^{^{\prime }}(x)\right) ^{2}}=\dfrac{20~x^{3}f_{c}(x)}{\left(
f_{c}^{^{\prime }}(x)\right) ^{2}}.
\end{equation*}%
If $N_{f_{c}}^{\prime }(x)=0$ we have $x=0$ or $f_{c}(x)=0$.

As the roots of $f_{c}(x)$\ are super-stable fixed points $(f^{\prime}(x) =
0)$ the only interesting critical point of $N_{f_{c}}$ is $0$ and we denote
it by $d_{2},$ so for the study of the iteration of $N_{f_{c}}$ we will
start at $x_{0}=d_{2}$.

\begin{figure}[th]
\vspace{5mm} \centerline{\epsfig{file=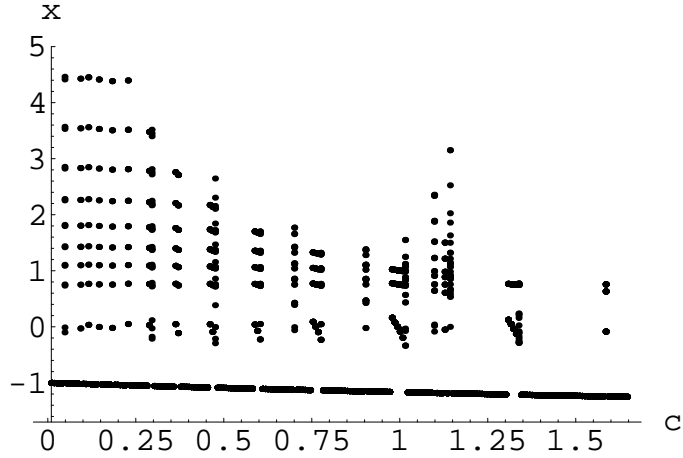,width=9cm}}
\caption{Bifurcation diagram of $N_{f_{c}}$ with $0<c<c_{0}$.}
\label{fig406}
\end{figure}

Let us now describe the numerical experiments which were performed in the $c$%
-parameter plane, see Figure \ref{fig406}.

To investigate this behavior further, we compute the bifurcation diagram for 
$N_{f_{c}}$ with $f_{c}(x)=x^{5}-c~x+1$, varying $c$ from $0$ to $c_{0}$.

\begin{remark}
In Figure \ref{fig406}\ we see that there is a sequence of \textquotedblleft
windows\textquotedblright\ where $N_{f_{c}}(d_{2})$ converges to a stable
periodic orbit with period $n$ $(n\in \mathbb{N)}$, intercalated with
intervals where the critical point $d_{2}$ converges to the fixed point $%
d_{0}$.
\end{remark}

Until now we have studied the case with only two coefficients in the quintic
polynomials (in such a case we have at most three roots).

Now we refer to the other two cases, i.e.:

\begin{enumerate}
\item the quintic polynomial with four distinct real roots, one of them
double;

\item the quintic polynomial with five distinct real roots.
\end{enumerate}

Now we recall the follow result:

\begin{theorem}[Barna \protect\cite{Barna53}]
If $f$ is a real polynomial having all real roots and at least four distinct
ones, then the set of initial values for which Newton's method does not
yield to a root of $f$ is homeomorphic to a Cantor set. The set of
exceptional initial values $J(f)$ is of Lebesgue measure zero.
\end{theorem}

\begin{proof}
For a proof we refer to Hurley and Martin \cite{Hurley84}. They all give
modern proof of Barna%
\'{}%
s result \cite{Barna53}. The underlying idea is to show that the set $J(f)$
arises in a manner which is very similar to the usual Cantor set
construction. Wong proves this result using symbolic dynamics \cite{Wong84}.
\end{proof}

In the next sections we concentrate in the most interesting case $%
f_{c}(x)=x^{5}-c~x+1$ for $c\in ]0,c_{0}[.$

\section{Symbolic dynamics}

Kneading theory is an appropriate tool to classify topologically the
dynamics of maps.

First we introduce the symbolic dynamics for the map $N_{f_{c}}$ where $%
f_{c}(x)=x^{5}-c~x+1$, and $0<c<c_{0}.$

\begin{figure}[th]
\vspace{5mm} \centerline{\epsfig{file=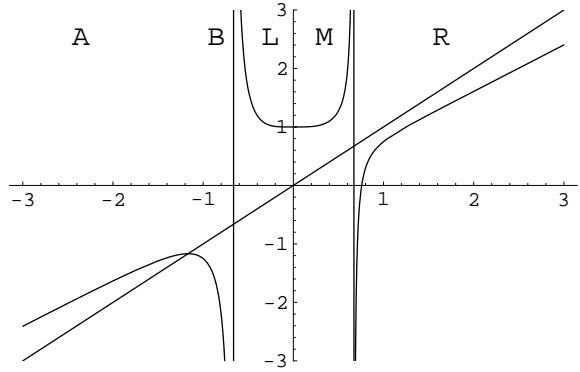,width=8cm}}
\caption{Symbolic dynamic for $N_{f_{c}}$ with $f_{c}=x^{5}-cx+1$ and $%
0<c<c_{0}$.}
\label{fig407}
\end{figure}

We consider the alphabet $\mathcal{A}=\{A,B,L,C,M,R\}$, and the set $\Omega =%
\mathcal{A}^{\mathbb{N}_{0}}$ of symbolic sequences on the elements of $%
\mathcal{A}$. Now we introduce the map 
\begin{equation*}
i_{c}:\mathbb{R}\backslash \underset{n\in \mathbb{N}_{0}}{\tbigcup }%
N_{f_{c}}^{-n}\left( \left\{ d_{1},d_{3}\right\} \right) \rightarrow \Omega
\end{equation*}%
defined by

\begin{equation*}
i_{c}(x)_{m}=\left\{ 
\begin{array}{l}
A\text{ \ if }N_{f_{c}}^{m}(x)<d_{0} \\ 
\\ 
B\text{ \ if }d_{0}<N_{f_{c}}^{m}(x)<d_{1} \\ 
\\ 
L\text{ \ if }d_{1}<N_{f_{c}}^{m}(x)<d_{2} \\ 
\\ 
C\text{ \ if }N_{f_{c}}^{m}(x)=d_{2} \\ 
\\ 
M\text{ \ if }d_{2}<N_{f_{c}}^{m}(x)<d_{3} \\ 
\\ 
R\text{ \ if }N_{f_{c}}^{m}(x)>d_{3}%
\end{array}%
\right.
\end{equation*}%
as we can see in Figure \ref{fig407}.

If we now consider the shift operator $\sigma :\Omega \rightarrow \Omega $, $%
\sigma (X_{0}X_{1}X_{2}...)=X_{1}X_{2}X_{3}...$\thinspace , we have the
commutative diagram%
\begin{equation*}
\begin{array}{ccccc}
&  & N_{f_{c}}\text{ } &  &  \\ 
& \Lambda & \longrightarrow & \Lambda &  \\ 
i_{c} & \downarrow &  & \downarrow & i_{c} \\ 
& \Sigma & \longrightarrow & \Sigma &  \\ 
&  & \sigma &  & 
\end{array}%
,
\end{equation*}%
where 
\begin{equation*}
\Lambda =\mathbb{R}\backslash \tbigcup\limits_{n\in \mathbb{N}%
}N_{f_{c}}^{-n}\left( \left\{ d_{1},d_{3}\right\} \right) .
\end{equation*}

We introduce in $\Omega $ an order, induced lexicographically by the order
in $\mathbb{R}$, with parity introduced by the subintervals where the
function is decreasing, $A<B<L<C<M<R$ when it is even and $-R<-M<C<-L<-B<-A $
when it is odd.

\begin{definition}
We say that $X\prec Y$ for $X,Y\in \Omega $, iff: 
\[
\exists _{k}:X_{i}=Y_{i},\forall _{0\leq i<k}\text{ and }%
(-1)^{n_{BL}(X_{1}...X_{k-1})}X_{k}<(-1)^{n_{BL}(X_{1}...X_{k-1})}Y_{k}
\]%
where $n_{BL}(X_{1}...X_{k-1})$ is equal to the number of times that the
symbols $B$ or $L$ appear in $X_{1}...X_{k-1}$.
\end{definition}

\begin{example}
$MRRM...\prec MRRR...$ and $RLRA...\succ RLRR...$
\end{example}

\begin{proposition}
Let $x,y\in \Lambda $. Then i) $x<y\Longrightarrow i_{c}(x)\preceq i_{c}(y)\ 
$and ii) $i_{c}(x)\prec i_{c}(y)\Longrightarrow x<y.$
\end{proposition}

\begin{proof}
It is sufficient to adapt the proof given to this end in Milnor and Thurston 
\cite{MT86}.
\end{proof}

We define the kneading sequence of the orbit of the critical point $x=d_{2}$
by

\begin{eqnarray*}
i &:&J\rightarrow \Omega \\
c &\mapsto &\sigma (i_{c}(d_{2})).
\end{eqnarray*}%
with 
\begin{equation*}
J=\{c:c\in \left] 0,c_{0}\right[ \text{ and }c\text{ is such that }\underset{%
n\in \mathbb{N}_{0}}{\cup }N_{f_{c}}^{n}(d_{2})\cap \left\{
d_{1},d_{3}\right\} =\emptyset \}.
\end{equation*}

Also we define the kneading sequences of the orbits of the discontinuous
point $d_{1}$ (respectively $d_{3}$) by $\sigma (i_{c}(d_{1}))$
(respectively $\sigma (i_{c}(d_{3}))$). We denote by $(U,X,Y,Z)$ the
kneading data $(\sigma (i_{c}(d_{0}))$, $\sigma (i_{c}(d_{1}))$, $\sigma
(i_{c}(d_{2}))$, $\sigma (i_{c}(d_{3})))$.

Now we characterize the admissible sequences looking at the typical graph of 
$N_{f_{c}}^{n}$ (see Figure \ref{fig407}). We get the following transition
matrix $T$ where rows and columns are labeled by the elements of $\mathcal{A}
$. 
\begin{equation}
T=\left[ 
\begin{array}{ccccc}
1 & 0 & 0 & 0 & 0 \\ 
1 & 0 & 0 & 0 & 0 \\ 
0 & 0 & 0 & 1 & 1 \\ 
0 & 0 & 0 & 1 & 1 \\ 
1 & 1 & 1 & 1 & 1%
\end{array}%
\right] .  \label{eq409}
\end{equation}

Then as the critical point $d_{2}$ is a local minimum, we get the following
admissibility 
\begin{equation}
\left\{ 
\begin{array}{l}
\sigma ^{i}(Y)_{1}=A\Rightarrow \sigma ^{i+1}(Y)=A^{\infty } \\ 
\\ 
\sigma ^{i}(Y)_{1}=B\Rightarrow \sigma ^{i+1}(Y)=A^{\infty } \\ 
\\ 
\sigma ^{i}(Y)_{1}=L\Rightarrow \sigma ^{i+1}(Y)\geq Y \\ 
\\ 
\sigma ^{i}(Y)_{1}=M\Rightarrow \sigma ^{i+1}(Y)\geq Y%
\end{array}%
.\right.  \label{eq410}
\end{equation}

Let the set $\Omega ^{+}=\{Y\in \mathcal{A}^{\mathbb{N}}:Y$ verifies $%
T_{Y_{i},Y_{i+1}}=1$ and (\ref{eq410})$\}$. We call $\Omega ^{+}$ the set of
admissible sequences.

\begin{example}
To see the admissibility we must pay attention to the fact that the critical
point $d_{2}$ is a local minimum and, in such case, if we have $\sigma
^{i}(Y)_{1}=L$ or $\sigma ^{i}(Y)_{1}=M$ (where $Y$ is a periodic sequence
with period $n$ of the critical point $d_{2}$, $1\leq i<n$) then we must
have $\sigma ^{i}(Y)>Y$. So the sequence $(RLRC)^{\infty }$ is admissible --
its occurrence can be seen in Figure \ref{fig406} near $c=1.3346...$ --
while the sequences $(LMAC)^{\infty }$ and $(RMRC)^{\infty }$ are not
admissible for the same reason.
\end{example}

\begin{figure}[th]
\vspace{5mm} \centerline{\epsfig{file=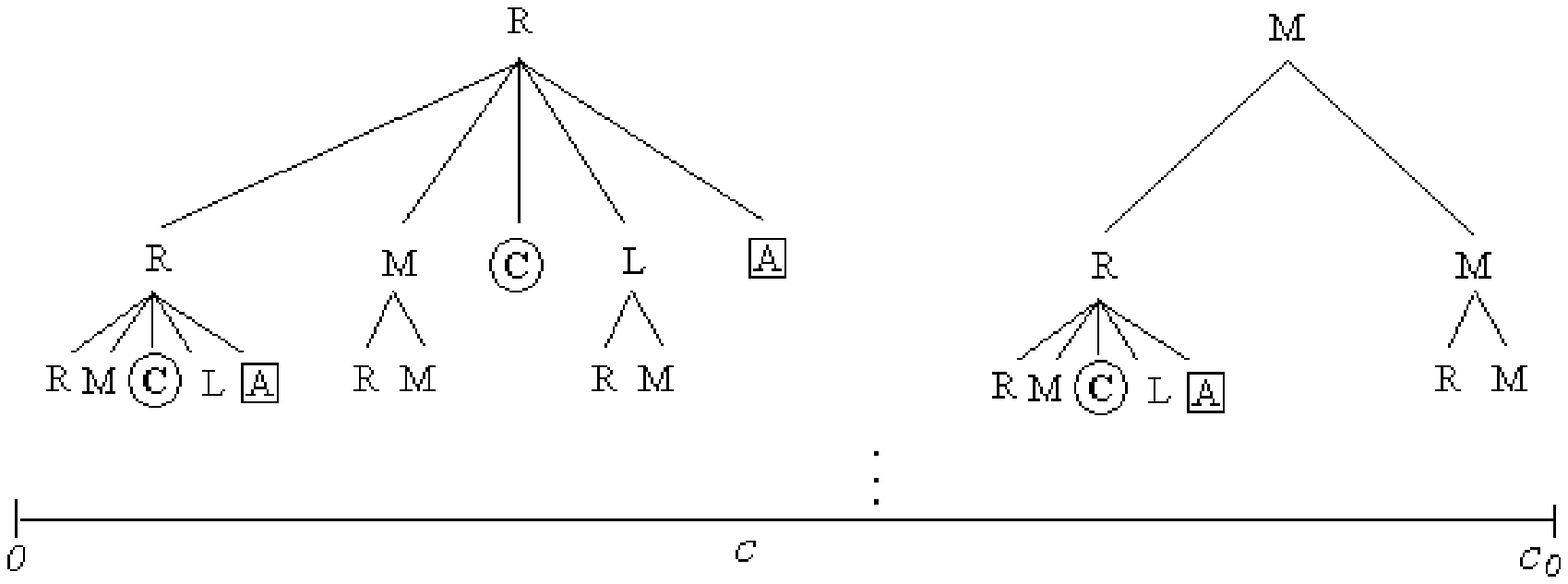,width=12cm}}
\caption{Tree $\mathcal{T}_{Y}$ of symbolic sequences corresponding to $%
d_{2} $ in the inverse order. The kneading sequences are marked with a
circle when the last symbol is $C$ and a square when the last symbol is $A$.}
\end{figure}

In a similar way to what is made by Milnor and Thurston in \cite{MT86} we
define the kneading increments associated to the kneading data by%
\begin{equation*}
\nu _{d_{i}}=\theta (d_{i}^{+})-\theta (d_{i}^{-})\text{ with }i=0,1,2,3
\end{equation*}%
where $\theta (x)$ is the invariant coordinate of each symbolic sequence
associated to the itinerary of each point $d_{i}$, see \cite{MT86}.

Using this we define the kneading matrix $N(t)$ and the kneading determinant 
\begin{eqnarray*}
D(t) &=&\frac{(-1)^{i+1}D_{i}(t)}{(1-\varepsilon _{i}t)} \\
&=&\frac{d_{Y}(t)}{(1-t)(1-t^{k})}
\end{eqnarray*}%
where $D_{i}(t)$ is the determinant of $N(t)$ without the column $i$ and the
cyclotomic polynomials in the denominator correspond to the stable periodic
orbits of $d_{0}$ and $d_{2}$, see \cite{MT86}.

\begin{example}
We exemplify the kneading increment for the sequence $RLRC$. We have 
\[
\theta (d_{0}^{+})=B-A t-A t^{2}-...=B-A t(1+t+t^{2}+...)=B-\frac{A t}{1-t}
\]%
and 
\[
\theta (d_{0}^{-})=A+At+At^{2}+...=A(1+t+t^{2}+...)=\frac{A}{1-t},
\]%
consequently%
\[
\nu _{d_{0}}=\theta (d_{0}^{+})-\theta (d_{0}^{-})=B-A\frac{1+t}{1-t}.
\]

Similarly, we get%
\[
\theta (d_{1}^{+})=L-Rt-Rt^{2}-...=L-(1+t+t^{2}+t^{3}+...)~Rt=L-\frac{Rt}{1-t%
}
\]%
and 
\[
\theta (d_{1}^{-})=B-At-At^{2}-...=B-(1+t+t^{2}+t^{3}+...)~At=B-\frac{At}{1-t%
},
\]%
consequently 
\[
\nu _{d_{1}}=\theta (d_{1}^{+})-\theta (d_{1}^{-})=L-B+\frac{(A-R)t}{1-t}.
\]%
Similarly, we get%
\begin{eqnarray*}
\theta (d_{2}^{+}) &=&M+Rt+Lt^{2}-Rt^{3}-Lt^{4}+... \\
&=&M+(1+t^{4}+t^{8}+...)(Rt+Lt^{2}-Rt^{3}-Lt^{4}) \\
&=&M+\frac{Rt+Lt^{2}-Rt^{3}-Lt^{4}}{1-t^{4}}
\end{eqnarray*}%
and%
\begin{eqnarray*}
\theta (d_{2}^{-}) &=&L-Rt-Lt^{2}+Rt^{3}+Lt^{4}-... \\
&=&L-(1+t^{4}+t^{8}+...)(Rt+Lt^{2}-Rt^{3}-Lt^{4}) \\
&=&L-\frac{Rt+Lt^{2}-Rt^{3}-Lt^{4}}{1-t^{4}},
\end{eqnarray*}%
consequently 
\[
\nu _{d_{2}}=\theta (d_{2}^{+})-\theta (d_{2}^{-})=M-L+\frac{%
2~(Rt+Lt^{2}-Rt^{3}-Lt^{4})}{1-t^{4}}.
\]%
Similarly, we get 
\[
\theta (d_{3}^{+})=R+At+At^{2}+...=R+(1+t+t^{2}+t^{3}+...)~At=R+\frac{At}{1-t%
}
\]%
and%
\[
\theta (d_{3}^{-})=M+Rt+Rt^{2}+...=M+(1+t+t^{2}+t^{3}+...)~Rt=M+\frac{Rt}{1-t%
},
\]%
consequently 
\[
\nu _{d_{3}}=\theta (d_{3}^{+})-\theta (d_{3}^{-})=R-M+\frac{(A-R)t}{1-t}.
\]%
So we have the kneading matrix 
\[
N(t)=\left[ 
\begin{array}{ccccc}
-\frac{1+t}{1-t} & 1 & 0 & 0 & 0 \\ 
-\frac{t}{1-t} & -1 & 1 & 0 & -\frac{t}{1-t} \\ 
0 & 0 & -1+\frac{2t^{2}-2t^{4}}{1-t^{4}} & 1 & \frac{2t-2t^{3}}{1-t^{4}} \\ 
\frac{t}{1-t} & 0 & 0 & -1 & 1-\frac{t}{t-1}%
\end{array}%
\right] .
\]%
With $i=2$ we have $\varepsilon _{2}=-1$ (because $N_{f_{c}}^{^{\prime
}}(x)|_{[d_{0},d_{1}]}<0$)%
\begin{eqnarray*}
D(t) &=&\frac{(-1)~D_{2}(t)}{1+t} \\
&=&\frac{(1+t)(1-t-t^{2}-t^{3})}{(1-t)(1-t^{4})}.
\end{eqnarray*}
\end{example}

Next we denote by $d_{Y}(t)$ the numerator of $D(t)$ given by $%
D(t)~(1-t)~(1-t^{k})$, where $k$ is the period of the critical point $d_{2}$
and $Y$ is the kneading sequence associated to $d_{2}$. Each kneading data
determines a kneading determinant but the most significant factor of the
numerator is determined by the kneading sequence $Y$.

It is easy to see the following result

\begin{proposition}
To the set $\Omega ^{+}$\ of the ordered kneading sequences can be
associated the tree $\mathcal{T}_{Y}$, where in each $k$-level of the tree
are localized kneading sequence of $~k$ length.

It is easy to see the following result
\end{proposition}

\begin{corollary}
To $\mathcal{T}_{Y}$ we associate a tree $\mathcal{T}_{d_{Y}(t)}$ of the
numerators of kneading determinant.
\end{corollary}

To proof this corollary we need the following lemma.

\begin{lemma}
Let $Y$ be an admissible periodic sequence corresponding to orbit of the
critical point $d_{2}$ of period $k$ whose the numerator of the kneading
determinant $d_{Y}(t),$ now we designate by $d_{k}(t).$ Then $d_{k}(t)$ has
degree $n=k$ and the polynomials correspondent to the periodic sequences of
period $k+1$ ($k+1$ - level of the tree) follow the rule of construction:%
\[
\begin{array}{c}
(1-t)d_{k}(t)=\underset{p(t)}{\underbrace{%
1-t+a_{2}t^{2}+a_{3}t^{3}+...+a_{k}t^{k-1}}}\underset{q(t)}{\underbrace{%
-\delta t^{k}-\delta t^{k+1}}}=p(t)+q(t) \\ 
\begin{array}{cccc}
\text{ \ \ \ \ \ \ \ \ \ \ \ \ \ \ \ \ \ }\swarrow & \downarrow & \downarrow
& \searrow \text{ \ \ \ } \\ 
\text{ \ \ \ \ \ \ \ \ \ \ \ }A & L & M & R \\ 
(1-t)d_{k+1}(t)=p(t)-2\delta t^{k};\text{ } & p(t)+2\delta t^{k};\text{ } & 
p(t)-2\delta t^{k+1};\text{ } & p(t)+tq(t)%
\end{array}%
\end{array}%
\]%
\textit{with} $a_{k}\in \{-2,0,2\}$ \textit{and} $\delta =(-1)^{n_{L}}$
where $n_{L}$ is equal to the number of times that the symbol $L$ appears in
the symbolic sequence $Y.$
\end{lemma}

\begin{proof}
Computing the kneading determinants of the sequences in each level $k$ of
the set $\Omega ^{+}$\ of the ordered kneading sequences and analyzing the
passage from level $k$ to level $k+1$ and using the induction it goes out
the rule of the indicated construction.
\end{proof}

\begin{remark}
If we want to see the symbolic dynamics in the tree for the converging
points in the Newton's method, it is of ramifications ending with $A^{\infty
}$. In this case $d_{Y}(t)$ stays constant after reaching the first symbol $%
A $ and we note when it reaches the symbol $A$ it is $A$ for ever.
\end{remark}

\begin{theorem}
\label{main}Let $P$ and $Q$ be kneading sequences in $\Omega ^{+}$ with the
lexicographic order $\prec $. If $P\prec Q$ then $c_{P}>c_{Q}$, where $c_{P}$
(respectively $c_{Q}$) is the parameter value corresponding to the kneading
sequence $P$ (respectively $Q$). If $c_{1}>c_{2}$ then there are $P,Q\in
\Omega ^{+}$ with $P\preceq Q$, where $P$ (respectively $Q$) is the kneading
sequence realized by the parameter value $c_{1}$ (respectively $c_{2}$).
\end{theorem}

\begin{proof}
It is sufficient to extend the Tsujii results on the quadratic map to the
Newton map. Then we prove that the kneading sequence $P$ associated to the
orbit of the critical point $d_{2}$ (minimum) is monotone decreasing with
respect to parameter $c.$
\end{proof}

\section{Topological entropy}

In the known paper by Misiurewicz and Szlenk \cite{Misiurewicz80} the
topological entropy is determined by 
\begin{equation}
h_{top}(N_{f_{c}})=\log s(N_{f_{c}}),  \label{eq406}
\end{equation}%
where $s(N_{f_{c}})=\underset{k\rightarrow \infty }{\lim }(L_{k})^{1/k}$, $%
L_{k}$ is the number of laps of $N_{f_{c}}^{k}$, i.e., the numbers of
sub-intervals where $N_{f_{c}}^{k}(x)$ is monotone (see also \cite{MT86}).

When the orbit of the critical point $d_{2}$ of $N_{f_{c}}(x)$ is periodic
we have a Markov partition which is determined by the itineraries of the
critical point. Once we have the Markov partition, a subshift of finite type
is determined by the transition matrix. Given a Markov partition $\mathcal{P}%
=\{I_{j}\}_{j=1}^{m}$, the transition matrix $\mathcal{M}=(a_{ij})$ of the
type $(n\times n)$ is defined by 
\begin{equation*}
a_{ij}=\left\{ 
\begin{array}{cc}
1 & \text{if ~}~int(N_{f_{c}}(I_{i})\cap I_{j})\neq \emptyset \\ 
0 & \text{if ~}~int(N_{f_{c}}(I_{i})\cap I_{j})=\emptyset%
\end{array}%
\right. .
\end{equation*}

Like in \cite{LSR93} the topological entropy $h_{top}(N_{f_{c}})$ is
obtained from the smallest real root $t^{\ast }$ of $d_{\mathcal{M}}(t)$, $%
t^{\ast }\in \lbrack \sqrt{2}-1,1]$, where $d_{\mathcal{M}}(t)=\det (I-t~%
\mathcal{M})$\ is the characteristic polynomial of the transition matrix $%
\mathcal{M}$. The value $t^{\ast }=\sqrt{2}-1$ corresponds to $c=c_{0}$,
which occurs to the kneading sequence $M^{\infty }$.

Also we can compute $h_{top}(N_{f_{c}})=\log 1/t^{\ast }$, where $t^{\ast }$
is the minimal solution of $D(t^{\ast })=0$, with 
\begin{equation*}
D(t)=\frac{d_{Y}(t)}{(1-t)~(1-t^{k})}=\frac{d_{\mathcal{M}}(t)}{%
(1-t)~(1-t^{k})~(1-t)^{2}}
\end{equation*}%
the kneading determinant, see \cite{MT86}, \cite{LSR93}.

\begin{remark}
The tree $\mathcal{T}_{d_{Y}(t)}$ shows the relation between the symbolic
sequences, namely the characteristic polynomials of the transition matrix $%
\mathcal{M}$,\ and the topological entropy for each parameter $c$ of the
family of Newton map $N_{f_{c}}$\ for the quintic $f_{c}(x)=x^{5}-c$~$x+1$,
with $c$ between $0$ and $c_{0}$.
\end{remark}

\begin{figure}[th]
\vspace{5mm} \centerline{\epsfig{file=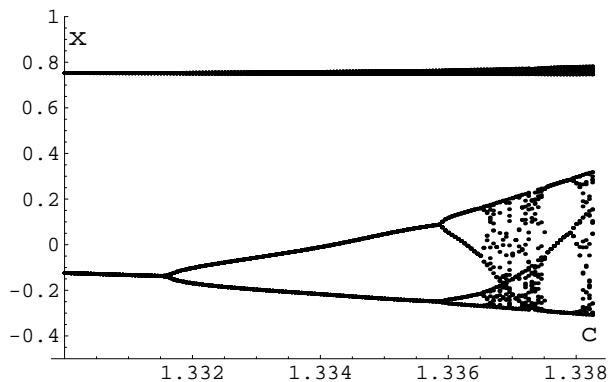,width=8cm}}
\caption{Bifurcation diagram near $c=1.3342$ to detect period duplication.}
\label{fig408}
\end{figure}

If we magnify the bifurcation diagram near $c=1.3342...$ we can see in
Figure \ref{fig408} the cascade beginning with period two on left and period
duplication.

We can illustrated the computation of the topological entropy for various
values of $c$.

\begin{example}
For $c=1.3342...$ we have period $4$ with the symbolic sequence $RLRC$%
\textit{. T}he Markov partition is in the Figure \ref{fig409}.

\begin{figure}[th]
\vspace{5mm} \centerline{\epsfig{file=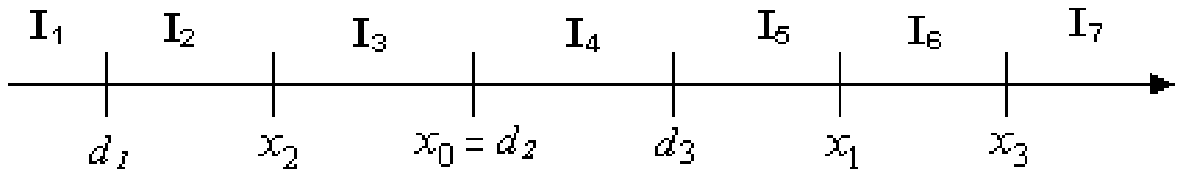,width=9cm}}
\caption{Markov partition for the sequence RLRC.}
\label{fig409}
\end{figure}

In this case with $c=1.3342...$ the transition matrix is 
\[
\mathcal{M}_{RLRC}=\left( 
\begin{array}{ccccccc}
1 & 0 & 0 & 0 & 0 & 0 & 0 \\ 
0 & 0 & 0 & 0 & 0 & 0 & 1 \\ 
0 & 0 & 0 & 0 & 0 & 1 & 0 \\ 
0 & 0 & 0 & 0 & 0 & 1 & 1 \\ 
1 & 1 & 0 & 0 & 0 & 0 & 0 \\ 
0 & 0 & 1 & 0 & 0 & 0 & 0 \\ 
0 & 0 & 0 & 1 & 1 & 1 & 1%
\end{array}%
\right)
\]%
and by formula (\ref{eq406}) we calculate $h_{top}(N_{f_{_{1.3342...}}})=%
\log 1.83929...$
\end{example}

\begin{figure}[th]
\vspace{5mm} \centerline{\epsfig{file=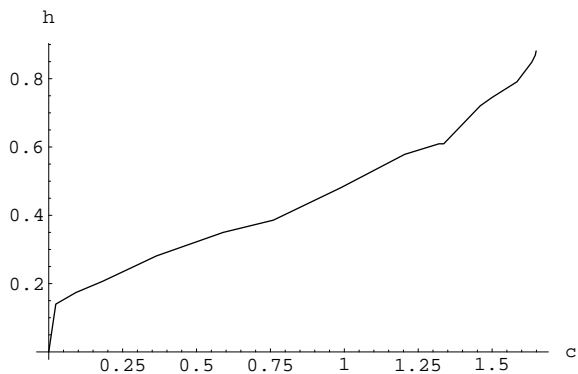,width=8cm}}
\caption{Entropy variation with respect to the parameter $c$ between $0$ and 
$c_{0}$.}
\label{fig410}
\end{figure}

We have computed more examples and we plot in Figure \ref{fig410}\ the graph
of the entropy for $c\in ]0,c_{0}[$.

\begin{corollary}
Let $N_{f_{c}}$ be the Newton map associated to the quintic map $%
f_{c}(x)=x^{5}-c$~$x+1$, with $c$ $\in ]0,c_{0}[$. The topological entropy
of $N_{f_{c}}$ is a non-decreasing function with respect to the parameter $c$%
. The topological entropy varies between $0$ and $\log (1+\sqrt{2})$.
\end{corollary}

\begin{proof}
This result is a consequence of the order and admissibility defined before
in the kneading sequence set $\Omega ^{+}$ and of Theorem \ref{main}. In the
extreme is the kneading sequence $M^{\infty },$ where $s(N_{f_{c}})=1+\sqrt{2%
}.$
\end{proof}

A next step for the future would be to study the types of bifurcation and
the Hausdorff dimension of the set of badly initial points, that is the set
of points not convergent to any root.

\section*{Acknowledgements}

The first author has been partially supported by D:G:E:S:I:C. Grant
PB-98-0374-C03-01. The second author is (partially) supported by CITMA
(Madeira Island - Portugal) trough the program POPRAM\ III and the third
author is (partially) supported by FCT (Portugal) trough the program
POCTI/FEDER.

\end{document}